\documentclass[12pt]{amsart}
\usepackage{graphicx}
\usepackage{amsfonts}
\usepackage{amssymb}
\usepackage{amsthm}
\usepackage[centertags]{amsmath}
\usepackage{newlfont}

\usepackage{color}
\usepackage{graphicx,color} 
\usepackage{amsmath, amssymb, graphics}

\def\r{\mathbb R}
\def\l{\mathbb L}
\def\c{\mathbb C}
\def\h{\mathbb H}
\def\t{\mathbf t}
\def\n{\mathbf n}
\def\b{\mathbf b}
\def\e{\mathbf E}

\newtheorem{theorem}{Theorem}[section]
\newtheorem{corollary}[theorem]{Corollary}

\newtheorem{remark}[theorem]{Remark}

\begin{document}

\title[New examples of maximal surfaces]{New examples of maximal surfaces in Lorentz-Minkowski space}
\author{Rafael L\'opez}
 \email{ rcamino@ugr.es}
 \address{Departamento de Geometr\'{\i}a y Topolog\'{\i}a\\ Instituto de Matem\'aticas (IEMath-GR)\\
 Universidad de Granada\\
 18071 Granada, Spain}
\thanks{The  first  author has been partially supported by  the MINECO/FEDER grant MTM2014-52368-P}

 \author{Seher Kaya}
\email {seher.Kaya@ankara.edu.tr}
   \address{Department of Mathematics, Ankara University\\ Ankara, Turkey}

\begin{abstract}
We use the  Bj\"{o}rling problem in   Lorentz-Minkowski space to obtain explicit parametrizations of maximal surfaces containing a circle and a helix.  We investigate the Weierstrass  representation of these surfaces.

\end{abstract}

\maketitle

\section{Introduction}
 In  three-dimensional Lorentz-Minkowski space  $\l^3$, the  Bj\"{o}rling problem consists of finding a maximal surface (spacelike surface with zero mean curvature) containing a given real analytic curve $\alpha:I\rightarrow\l^3$ and a prescribed unit normal vector field $V$ along $\alpha$. Here $V$ is a real analytic unit timelike vector field orthogonal to $\alpha'$ along $I$.   In \cite{acm} the authors solve affirmatively the problem. For this, we  view the interval $I$ as $I\times\{0\}\subset\mathbb{C}$ and  by analyticity, the functions $\alpha$ and $V$ have holomorphic extensions $\alpha(z)$ and $V(z)$ in a simply-connected domain $\Omega\subset\c$ that contains $I\times\{0\}$. Then the Bj\"{o}rling solution is the surface $X:\Omega\rightarrow\l^3$ whose parametrization is
\begin{equation}\label{bj}
X(u,v)=\mbox{Re}\left(\alpha(z)+i\int\limits_{u_0}^z V(w)\times \alpha'(w)\ dw \right)
\end{equation}
where $u_0\in I$ is fixed and $z\in\Omega$ (\cite{acm}).  The surface $X(u,v)$ is unique with the prescribed Bj\"{o}rling data $\alpha$ and $V$, and  $\alpha$ being the parameter curve $v=0$. We point out that  $\alpha$ is a spacelike curve because any curve contained in a spacelike surface must be spacelike. 

Although the Bj\"{o}rling formula \eqref{bj} only requires a complex integration, there are  few examples in the literature of explicit parametrizations of Bj\"{o}rling surfaces: see \cite{acm}. These surfaces are the Lorentzian counterpart of the known examples of minimal surfaces of Euclidean space \cite{gas}. The goal of this paper is to provide new explicit parametrizations of solutions of the Bj\"{o}rling problem. Our surfaces appear when the core curve  $\alpha$ is a circle or a helix of $\l^3$. The idea that lies behind these examples is the following.  Suppose that $\alpha$ is a Frenet curve of $\l^3$, that is, a curve  where the Frenet frame $\{\t,\n,\b\}$ is formed by an orthonormal basis of $\l^3$. Then the unit vector field $V(t)$ is a linear combination of $\n(t)$ and $\b(t)$ for all $t\in I$. Our examples appear when the curve $\alpha$ is a circle or a helix and the  linear combination that defines $V(t)$ is given by trigonometric polynomials on the variable $t$.  The simple but crucial step is that when the integrand in \eqref{bj}  is formed by trigonometric polynomials, then the integral in \eqref{bj} can be explicitly obtained. The above idea was developed in Euclidean space by Meeks and Weber in \cite{mw}. Indeed, when $\alpha$ is a circle, the Bj\"{o}rling surfaces are called bending helicoids \cite{mw} and when  $\alpha$ is a helix, the surfaces are called  helicoidal helicoids \cite{we}.

 In Lorentzian ambient space, we establish the same program and we follow here the same terminology  for maximal surfaces of $\l^3$.  However, there are some differences that need to point out.  First, there exist circles and helices in $\l^3$ that do not parametrize with trigonometric polynomial but as polynomial functions on the variable $t$. Exactly in $\l^3$ there are three types of circles depending on the causal character of the rotation axis, namely, the axis is timelike, spacelike or lightlike. The circles whose axis is lightlike are parameterized by polynomials on the variable. Thus our procedure applies   only  when the axis of the circle and the helix is timelike or spacelike.   A second difference is related with the vector field $V$. For each $t$, the vector $V(t)$ lies in  the  plane spanned by $\{\n(t),\b(t)\}$ whose induced metric is Lorentzian. Since $V(t)$ is a unit timelike vector, then  $V(t)$ belongs to the timelike unit circle determined by $\n(t)$ and $\b(t)$, which is a hyperbola for each plane $\mbox{span}\{\t(t)\}^\bot$. In particular, this implies that $V(t)$ does not take any direction in that plane, in contrast to what happens in the Euclidean case.

The explicit parametrizations of the new maximal surfaces are given in terms of trigonometric functions but they have cumbersome expressions. Exactly the parametrizations are:
\begin{enumerate}
\item The core curve $\alpha$ is a circle.
\begin{enumerate}
\item The axis is timelike: equation \eqref{ct}.
\item The axis is spacelike: equations \eqref{cs1} and \eqref{cs2}
\end{enumerate}
\item The core curve $\alpha$ is a helix.
\begin{enumerate}
\item The axis is timelike: equation \eqref{ht}.
\item The axis is spacelike: equations \eqref{hs11}, \eqref{hs12}, \eqref{hs21} and \eqref{hs22}.
\end{enumerate}
\end{enumerate}

\section{Maximal surfaces in $\l^3$}

 The Lorentz-Minkowski space   is the vector space $\r^3$ endowed with the Lorentzian metric $\langle,\rangle  =dx^2+dy^2-dz^2$, where $(x,y,z)$ are the canonical coordinates of $\r^3$. A smooth immersion $X:M\rightarrow\l^3$ of a (connected) surface $M$ is said to be  spacelike if the induced metric $ds^2$ on $M$ via $X$ is Riemannian. As a consequence of the spacelike condition, $M$ is an orientable surface because a Gauss map $N$   on $M$ is a unit timelike vector field and thus its range lies in one of the two connected components of the hyperbolic plane $\h^2=\{p\in\l^3: \langle p,p\rangle=-1\}$.    We assume that  $N:M\rightarrow\h^2_{+}$ where $\h_{+}^2=\{p\in\h^2: z(p)>0\}$.   A maximal surface is a spacelike immersion $X:M\rightarrow\l^3$ with zero mean curvature  everywhere and we also say that $M$ is a maximal surface if $X$ is understood.  
 
For a spacelike surface,  $(M,ds^2)$ is a Riemann oriented surface and thus $M$ is endowed with a natural conformal structure of a Riemann surface. If, in addition, $M$ is a maximal surface, we suppose that   $X:(M,ds^2)\rightarrow\l^3$ is a   conformal maximal immersion with  $X=X(z)$ and  $z=u+iv\in \c$ is a  conformal parameter. Because $X$ is maximal, the coordinate functions $X_k$ are harmonic and thus $dX_k$ are harmonic $1$-forms on $M$. Let $\phi_k=dX_k+i(\ast dX_k)$ the holomorphic $1$-forms that extend $dX_k$ with respect to the complex structure associated to $(M,ds^2)$. Then we have
$$\phi_1^2+\phi_2^2-\phi_3^2=0,\quad\quad ds^2=|\phi_1|^2+|\phi_2|^2-|\phi_3|^2.$$
 If $\phi=(\phi_1,\phi_2,\phi_3)$, the surface is then obtained by $X(p)=X(p_0)+\mbox{Re}(\int_\gamma\phi)$ for any curve $\gamma$ connecting $p_0$ and $p$. The integral does not depend on the curve $\gamma$ which is equivalent to  $\mbox{Re}\int_\gamma\phi=0$ for any closed curve $\gamma$ in $M$ and we say that $\phi_k$ have no real periods. Define
$$g=\frac{\phi_3}{\phi_1-i\phi_2},\quad\quad \omega=fdz=\phi_1-i\phi_2.$$
Then   $g$ is a meromorphic function on $M$ and $\omega$ a holomorphic $1$-form on $M$.  The pair $(g,\omega)$ is called the Weierstrass data of the surface and we can write  
$$\phi_1= \frac12 (1+g^2)\omega,\quad \phi_2=\frac12  i(1-g^2)\omega,\quad\phi_3=g\omega.$$
The  Gauss map of $X$   is
$$N=\frac{1}{1-|g|^2}\left(2 \mbox{Re} (g),2 \mbox{Im} (g),1+|g|^2\right),$$
and if $\pi$ is    the stereographic projection   from $\h^2_{+}$ into the unit disc $D\subset\c$ from the point $(0,0,-1)$, then $g=\pi\circ N$.

There is a reverse process to obtain maximal surfaces as follows. Consider $M$ a Riemann surface and let $X:M\rightarrow\l^3$ be a non-constant harmonic   map with $X=X(z)$. If $\phi_k =\partial X_k/\partial z$ ($1\leq k\leq 3$) satisfy $\phi_1^2+\phi_2^2-\phi_3^2=0$ and $|\phi_1|^2+|\phi_2|^2-|\phi_3|^2\not\equiv 0$, then $X:M\setminus \mathcal{R}\rightarrow\l^3$ is a maximal surface where $\mathcal{R}$ is the set of branch points of $X$, that is, points where $X$ is not an immersion or points where the induced metric from $\l^3$ is degenerate.  We say that $X:M\rightarrow\l^3$ is a generalized maximal surface \cite{er}.

Finally we point out that there exists a duality between maximal surfaces of $\l^3$ and minimal surfaces in Euclidean space $\e^3$ as a metric space in the sense that for each maximal surface of $\l^3$ it is possible to construct a minimal surface in $\e^3$ and {\it vice versa}. Indeed, if $(\phi_1,\phi_2,\phi_3)$ are the holomorphic $1$-forms of a maximal surface $M$ of $\l^3$, then
$(\psi_1,\psi_2,\psi_3)= (i\phi_1,i\phi_2,\phi_3)$  are the holomorphic $1$-forms of  a minimal surface $M^\sharp$ of $\e^3$. If $(g,\omega)$ is the Weierstrass data of $M$, then $(-ig,i\omega)$ is for $M^\sharp$. Reciprocally, if $\psi_k$ are the holomorphic  $1$-forms of a minimal surface of $\e^3$, then $(-i\psi_1,-i\psi_2,\psi_3)$ are the  holomorphic $1$-forms of a generalized maximal surface.


\section{Bj\"{o}rling surfaces based on a circle }

In this section we construct   Bj\"{o}rling surfaces  when the core curve $\alpha$ is a   circle, in particular, we extend in the Lorentzian ambient the bending helicoids of Meeks and Weber.   By a  \emph{circle}   in $\l^3$ we mean the orbit of a point under a uniparametric group of rotational  motions  when this orbit is not a straight line.   A circle is  determined by a point and the axis $L$ of the group of rotational motions.  

When the axis is timelike or spacelike axis, a spacelike circle parametrized by the arc-length is a Frenet curve and  the Frenet frame is $\t(t)=\alpha'(t)$, $\n(t)=\alpha''(t)/|\alpha''(t)|$ and $\b(t)=\t(t)\times\n(t)$. The causal character of $\n(t)$ and $\b(t)$) depends on the one of $\alpha''(t)$. Moreover,   $\b(t)$ is a constant vector field, namely, the direction of the  axis. If the axis of the circle is lightlike, then the curve  is not a Frenet curve because there is not a canonical Frenet frame given by an orthonormal basis. In this case, $\alpha''(t)$ is a lightlike vector and the Frenet frame is formed by $\alpha'(t)$ and two lightlike vectors: see subsection \ref{sl}.

We   distinguish the study according to the axis of the circle. After a change of coordinates we will assume that the axis $L$ of a circle is spanned by $(0,0,1)$, $(1,0,0)$ and $(1,0,1)$.
\subsection{The axis is timelike}

 After a homothety  of $\l^3$, a spacelike circle with timelike axis $(0,0,1)$  parametrizes as $\alpha(t)=(\cos(t),\sin(t),0)$, $t\in\r$. The normal and binormal vectors of $\alpha$ are
$$\n(t)=(-\cos(t),-\sin(t),0),\quad\quad \b(t)=(0,0,1).$$
 Consider $V(t)$   a timelike unit vector field along $\alpha$. Because $\{\n(t),\b(t)\}$ is an orthonormal basis of $\mbox{span}(\alpha'(t))^\bot$, then we can write 
\begin{equation}\label{vv}
V(t)=\sinh(\varphi(t)) \n(t)+\cosh(\varphi(t)) \b(t)
\end{equation}
  for some function $\varphi(t)$. The case that $\varphi$ is a constant function is known, exactly we have:
 \begin{enumerate}
 \item If $\varphi=0$, then the Bj\"{o}rling surface is a (horizontal) plane of equation $z=\mbox{constant}$. 
 \item If $\varphi(t)=a\in\r$, $a\not=0$, then the Bj\"{o}rling  surface is a rotational surface. Indeed, an immediate integration of \eqref{bj} gives
 $$X(u,v)=\left(\cos (u) (\cosh (a) \sinh (v)+  \cosh (v)), \sin (u)(\cosh (a) \sinh (v)+ \cosh (v)),v \sinh (a)\right).$$ 
In order to see that this surface is rotational, we recall that  the uniparametric group of rotations with axis $(0,0,1)$ is  
 $$G_t=\left\{\Psi^t(\theta)=
\left(
\begin{array}{ccc}
\cos(\theta)  & -\sin(\theta)&0\\
\sin(\theta) & \cos(\theta) & 0 \\
0 & 0&1 \\
\end{array}
\right):\theta\in\r\right\}.$$
Then we have $\Psi^t(\theta)\cdot X(u,v)=X(u+\theta,v)$ for all $\theta\in\r$, proving that $X(u,v)$ is rotational. This surface is called an elliptic catenoid (or a catenoid of first king according to   \cite{ko}).  
 \end{enumerate}
  
The bending helicoids of $\l^3$ are obtained when in the expression of $V(t)$ given in \eqref{vv}, we take $\varphi(t)=at$, $a>0$, and we solve the corresponding Bj\"{o}rling  solution.   Then
\begin{eqnarray*}
V(t)&=&\sinh(at)\n(t)+\cosh(at)\b(t)\\
&=&(-\sinh(at)\cos(t),-\sinh(at)\sin(t),\cosh(at)),\quad t\in\r.
\end{eqnarray*}
We extend analitically $V$ and $\alpha$  and denote their extensions as $V(z)$ and $\alpha(z)$.  Then
$$V(z)\times \alpha'(z)=(-\cosh (az) \cos (z),-\cosh (az) \sin (z),\sinh (az)).$$
Once integrated \eqref{bj}, the parametrization of the  Bj\"{o}rling  surface  is
\begin{eqnarray}\label{ct}
X(u,v)&=&\Big(\frac{a \cos (u) \sin (a v)\cosh (a u)\cosh (v)  }{a^2+1}-\frac{a \sin (u) \cos (a v)\sinh (a u)\sinh (v)  }{a^2+1}\nonumber\\
&& +\frac{\cos (u) \cos (a v)\cosh (a u) \sinh (v) }{a^2+1}+\frac{\sin (u) \sin (a v)\sinh (a u)\cosh (v)  }{a^2+1}\nonumber\\
&&+\cos (u) \cosh (v),\nonumber\\
&& \frac{a \cos (u)\cos (a v) \sinh (a u)\sinh (v)  }{a^2+1}+\frac{a \sin (u) \sin (a v)\cosh (a u)\cosh (v)  }{a^2+1}\nonumber\\
&&+\frac{\sin (u) \cos (a v) \cosh (a u)\sinh (v) }{a^2+1}-\frac{\cos (u) \sin (a v) \sinh (a u)\cosh (v) }{a^2+1}\nonumber\\
&&+\sin (u) \cosh (v),\nonumber\\
&& -\frac{ \sin (a v)\sinh (a u)}{a}\Big).
\end{eqnarray}
We compute the Weierstrass data. From \eqref{bj}, $\phi(z)=\alpha'(z)+i V(z)\times \alpha'(z)$. Hence we have
\begin{eqnarray*}
\phi_1&=&(-\sin (z)-i \cos (z) \cosh (a z))dz\\
&=& \left(-\frac{1}{4} i \left(e^{-i z}+e^{i z}\right) \left(e^{-a z}+e^{a z}\right)-\frac{1}{2} i \left(e^{-i z}-e^{i z}\right)\right)dz\\
\phi_2&=&(\cos (z)-i \sin (z) \cosh (a z))dz\\
&=&\left(\frac{1}{4} \left(e^{-i z}-e^{i z}\right) \left(e^{-a z}+e^{a z}\right)+\frac{1}{2} \left(e^{-i z}+e^{i z}\right)\right)dz\\
\phi_3&=&(i\sinh (a z) dz=\frac{1}{2} i (e^{-a z}- e^{a z}) dz
\end{eqnarray*}
and the Weierstrass representation is
$$
\omega =- i\frac{ \left(e^{a z}+1\right)^2}{2e^{(a+i) z}}dz,\quad
 g(z)=-\frac{e^{i z} \left(e^{a z}-1\right)}{e^{a z}+1}.$$

 \begin{remark} {\rm
\begin{enumerate}
\item Comparing with the bending helicoids in Euclidean space \cite{mw}, and for small values of $v$, our examples give {\it embedded} strips of maximal surfaces  only $\alpha$ is not a  full circle, that is, only when $\alpha$ is a piece of length less than $2\pi$ because $V(t)$ is not a periodic function. A way to get  a periodic vector field $V(t)$ is replacing the function $\varphi(t)$ in \eqref{vv}  by a periodic function, as for example, $\varphi(t)=\cos(t)$. However we have no trigonometric polynomial in \eqref{bj} to can be integrated.  
\item In \cite{mw}, Meeks and Weber obtained non-orientable surfaces by twisting one half the vector field $V(t)$ along the circle   (also \cite{mi}). This is not possible to do it  in the Lorentzian ambient because $V(t)\in\h_{+}^2$ and  $V(t)$ can not reverse its initial position.  This is expectable because any spacelike surface is orientable.
\end{enumerate}}
 \end{remark}

\begin{figure}[h]
\includegraphics[width=.4\textwidth]{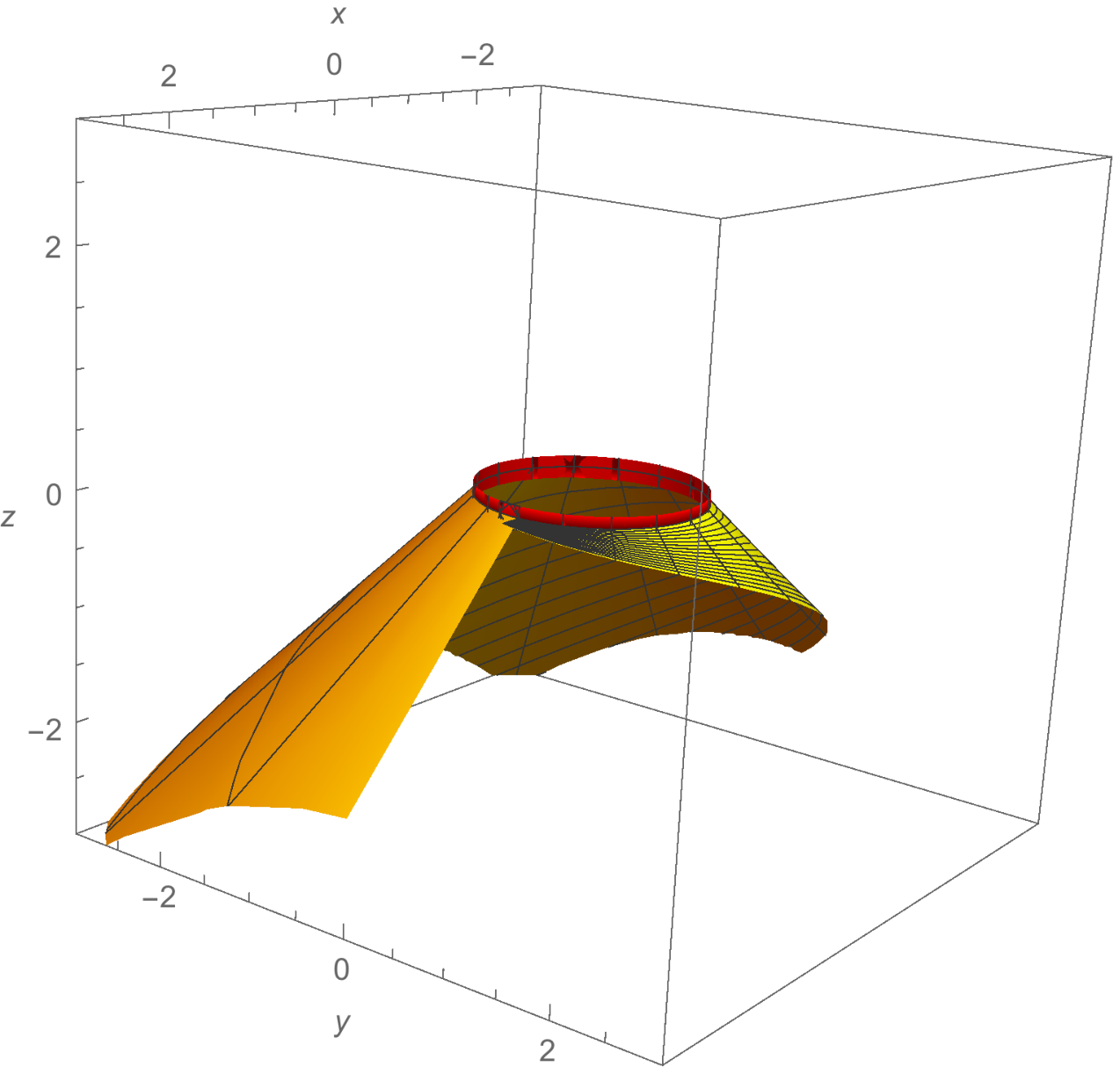}
\includegraphics[width=.4\textwidth]{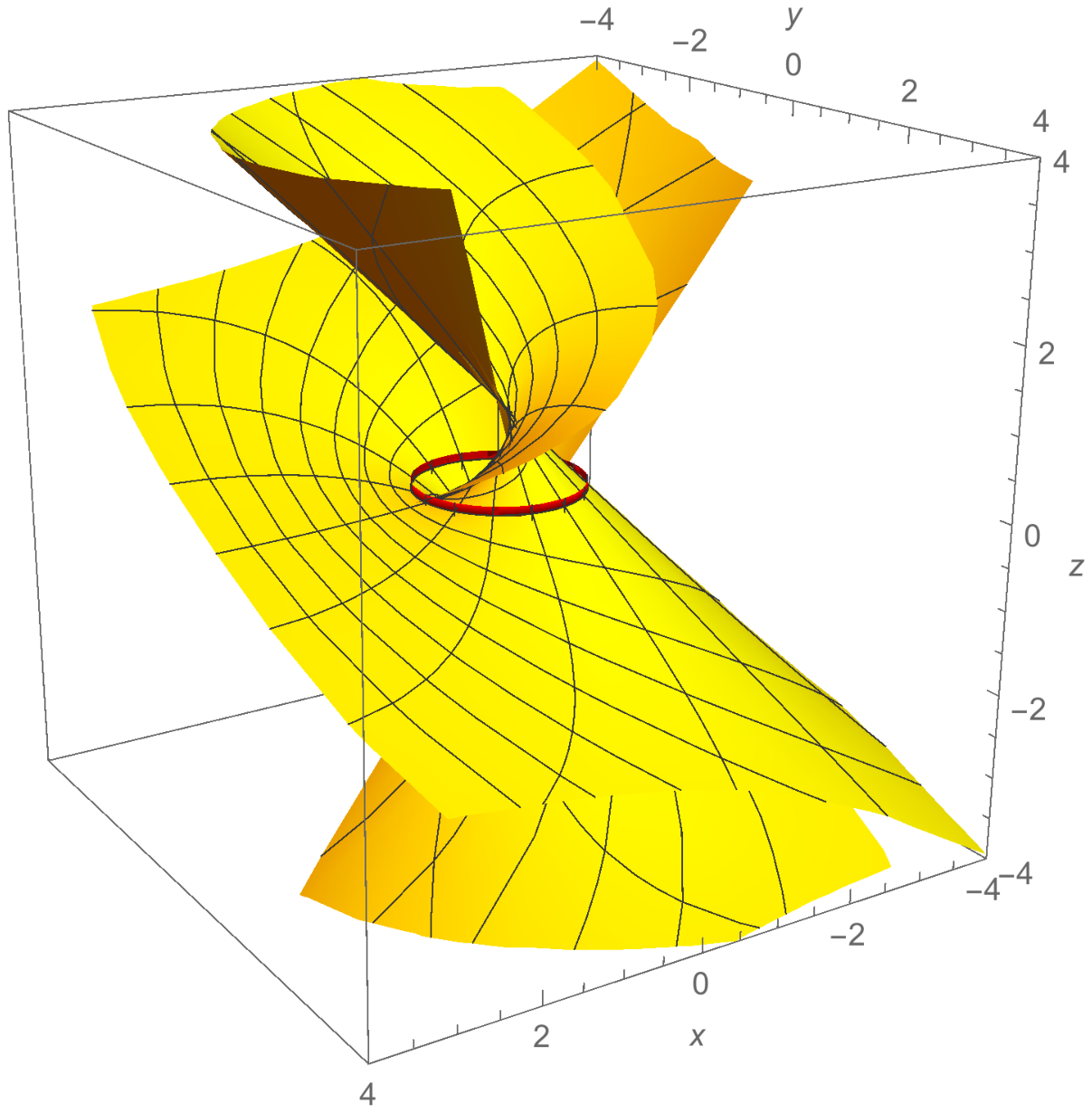}
\caption{A bending helicoid $X(u,v)$ where the core curve is a circle $\alpha$ with timelike axis. Here $a=1$. Left: a strip of $X(u,v)$ along $\alpha(u)=X(u,0)$. Right: a view of the full surface.}
\end{figure}
\subsection{The axis is spacelike}

Let $\alpha$ be a spacelike circle with spacelike axis $(1,0,0)$.  After a homothety, we parametrize the circle  as $\alpha(t)=(0,\sinh(t),\cosh(t))$, $t\in\r$. This curve is  a Frenet frame, namely, the normal and binormal vectors of $\alpha$ are
$$\n(t)=(0,\sinh(t),\cosh(t)),\quad\quad \b(t)=(1,0,0).$$
 Now a unit timelike vector field  $V(t)$ expresses as $V(t)=\cosh(\varphi(t)) \n(t)+\sinh(\varphi(t)) \b(t)$. The case   $\varphi(t)$ is constant gives a rotational surface. Indeed, if $\varphi(t)=a\in\r$,  a direct integration of the Bj\"{o}ling surface by using \eqref{bj} gives
 $$X(u,v)=\left(\begin{array}{l}
 v \cosh (a)\\
 \sinh (a) \sinh (u) \sin (v)+\sinh (u) \cos (v)\\
 \sinh (a) \cosh (u) \sin (v)+\cosh (u) \cos (v)
 \end{array}\right).$$
 On the other hand, the uniparametric group of rotations with axis $(1,0,0)$ is  
 $$G_s=\left\{\Psi^s(\theta)=
\left(
\begin{array}{ccc}
 1  & 0&0\\
0 & \cosh(\theta) & \sinh(\theta) \\
0 & \sinh(\theta)& \cosh(\theta) \\
\end{array}
\right):\theta\in\r\right\}.$$
It is easy to see that $\Psi^s(\theta)\cdot X(u,v)=X(u+\theta,v)$ for all $\theta\in\r$, proving that $X(u,v)$ is rotational.  This surface is called the hyperbolic catenoid. When $a=0$, the parametrized surface appears in \cite{ko} as  the catenoid of second kind.  
 
We now consider  bending helicoids when the core curve is the above circle and $\varphi(t)=at$, $a>0$. Then
\begin{eqnarray*}
V(t)&=&\cosh(at)\n(t)+\sinh(at)\b(t)\\
&=&(\sinh(at),\cosh(at)\sinh(t),\cosh(at)\cosh(t)),\quad t\in\r.
\end{eqnarray*}
Then extending analitically, we obtain
$$V(z)\times \alpha'(z)=(-\cosh (a z),-\sinh (z) \sinh (a z),-\cosh (z) \sinh (a z)).$$

The explicit integration of   \eqref{bj}  with $u_0=0$ depends if  $a\not=1$ or $a=1$. If $a\not=1$, the parametrization of the surface is
\begin{eqnarray}\label{cs1}
X(u,v)&=&\Big(\frac{\cosh (a u) \sin (a v)}{a},\nonumber\\
&&-\frac{ \sin (v) \cos (a v)\sinh (u)\sinh (a u) }{a^2-1}+\frac{a \sin (a v) \cos (v) \sinh (u)\sinh (a u) }{a^2-1}\nonumber\\
&& +\frac{a  \sin (v)  \cos (a v)\cosh (u)\cosh (a u)}{a^2-1}-\frac{ \cos (v)\sin (a v) \cosh (u)\cosh (a u) }{a^2-1}\nonumber\\
&&+\sinh (u) \cos (v),\nonumber\\
&&\frac{a  \cos (v)\sin (a v) \cosh (u)\sinh (a u) }{a^2-1}-\frac{ \sin (v)  \cos (a v)\cosh (u)\sinh (a u)}{a^2-1}\nonumber\\
&& +\frac{a  \sin (v)\cos (a v) \sinh (u)\cosh (a u) }{a^2-1}-\frac{\cos (v)\sin (a v)\sinh (u)  \cosh (a u) }{a^2-1}\nonumber\\
&&+\cos (v)\cosh (u) \Big).
\end{eqnarray}

When $a=1$,  the surface is
\begin{equation}\label{cs2}
X(u,v) = \left(\begin{array}{l}
\cosh (u) \sin (v)\\
\sinh (u) \cos (v)+\frac{1}{2} \sin (v) \cos (v)(\sinh ^2(u) + \cosh ^2(u))-\frac{v}{2}\\
\cosh (u) \cos (v) (1+\sinh (u) \sin (v) )
\end{array}\right).
\end{equation}

The  holomorphic $1$-forms $\phi_k$ are
\begin{eqnarray*}
\phi_1&=&-i \cosh (a z)dz=-\frac{1}{2} i e^{-a z} \left(e^{2 a z}+1\right)dz\\
\phi_2&=& \cosh (z)-i \sinh (z) \sinh (a z)dz\\
&=&\frac{1}{2} \left(e^{-z}+e^z\right)-\frac{1}{4} i \left(e^{z}-e^{-z}\right)\left(e^{a z}-e^{-a z}\right)dz\\
\phi_3&=&\sinh (z)-i \cosh (z) \sinh (a z)dz\\
&=&\frac{1}{2} \left(e^z-e^{-z}\right)-\frac{1}{4} i \left(e^{z}+e^{-z}\right) \left(e^{a z}-e^{-a z}\right)dz
\end{eqnarray*}
and the Weierstrass data are
$$\omega= -\frac{\left(e^{(a+1) z}+i e^{a z}+i e^z+1\right)^2}{4 e^{(a+1)z} }dz,\quad
g(z)=\frac{i e^{(a+1)z}+e^{a z}-e^z-i}{e^{(a+1)z}+i e^{a z}+i e^z+1}.$$

In contrast to the case of timelike axis,   the expressions of $\phi_k$ are only given  in terms of the exponential function $e^z$. This   allows to do the change of variable $e^z\rightarrow z$ (and $dz\rightarrow dz/z$), obtaining
\begin{eqnarray}\label{bhspacelike}
\phi_1&=&-  i\frac{z^{2 a}+1}{2z^{a+1}}dz\nonumber\\
\phi_2&=& \frac{-i z^{2 a+2}+i z^{2 a}+2 z^{a+2}+2 z^a+i z^2-i}{4z^{a+2}}dz\\
\phi_3&=&  \frac{-i z^{2 a+2}-i z^{2 a}+2 z^{a+2}-2 z^a+i z^2+i}{4z^{a+2}}dz.\nonumber
\end{eqnarray}
The Weierstrass representations is
$$\omega=-\frac{\left(z^{a+1 }+i z^a+i z+1\right)^2}{4 z^{a+1} }dz\quad\quad g(z)=i\frac{z^{a+1}-iz^a+iz-1}{z^{a+1}+iz^a+i z+1}.$$
The  holomorphic $1$-form $\omega$ is defined on $M=\c-\{0\}$. When $a=n\in\mathbb{N}$,   it is clear from \eqref{bhspacelike} that $\phi_k$ have not real periods if $n\not=1$ and when $n=1$, it is only $\phi_2$ that has are real periods.

 We end this section relating the maximal surfaces obtained in this subsection with    minimal surfaces of Euclidean space. Denote by $M_{n,s}$ the bending helicoid for spacelike axis obtained from  \eqref{bhspacelike} for $a=n\in\mathbb{N}$. Then the holomorphic $1$-forms of  $M_{n,s}^\sharp\subset\e^3$ are   
\begin{eqnarray*}
\psi_1&=&\frac{z^{2 n}+1}{2z^{n+1}}dz\\
\psi_2&=& \frac{z^{2 n+2}- z^{2 n}+2i z^{n+2}+2i z^n- z^2+1}{4z^{n+2}}dz\\
\psi_3&=&  \frac{-i z^{2 n+2}-i z^{2 n}+2 z^{n+2}-2 z^n+i z^2+i}{4z^{n+2}}dz.
\end{eqnarray*}
The Weierstrass data of  $M_{n,s}^\sharp$ are 
$$\omega=-i\frac{\left(z^{n+1 }+i z^n+i z+1\right)^2}{4 z^{n+1} }dz,\quad\quad g(z)=\frac{z^{n+1}-iz^n+iz-1}{z^{n+1}+iz^n+i z+1}.$$
 Then the $1$-forms $\psi_k$ are defined  on $\c-\{0\}$ and they have  no real periods. It is immediate that $M_{n,s}^\sharp$ is complete and with total finite curvature $-4\pi (n+1)$. The surface   $M_{n,s}^\sharp$ has two ends corresponding to $z=0$ and $z=\infty$. We focus in the case of least total curvature, that is, for  $n=1$, where the total curvature of $M_{1,s}^\sharp$  is $-8\pi$ and
$$\psi=\left(\frac{z^2+1}{2z^2},\frac{z^4+2i z^{3} - 2z^2+2iz+1}{4z^3},  \frac{-iz^4+2z^3-2z+i}{4z^3}\right)dz.$$
$$\omega=-i\frac{\left(z^{2 }+2i z +1\right)^2}{4 z^{2} }dz\quad\quad g(z)=\frac{z^{2}- 1}{z^{2 }+2i z +1}.$$
The orders of the ends $z=0$ and $z=\infty$ are $\nu_0=\nu_\infty=3$. This implies that the ends are asymptotic to an Enneper end, in particular, they are not embedded.  Minimal surfaces with Enneper ends are well known (\cite{cg,es,fs,ka,sa,th}) and when the genus  of the surface is $0$, these surfaces were studied by Karcher in \cite{ka}. When $n=1$, the surface $M_{1,s}^\sharp$ appears  in the classification  of the complete minimal surfaces with total curvature $-8\pi$ done by L\'opez, see \cite[Table 1]{lo}.

\begin{corollary} For $n=1$, the minimal surface $M_{1,s}^\sharp\subset\e^3$ corresponding to the bending helicoid $M_{1,s}$ of spacelike axis has  with total curvature $-8\pi$ and  two Enneper ends of order $3$.
\end{corollary}

\begin{figure}[h]
\includegraphics[width=.4\textwidth]{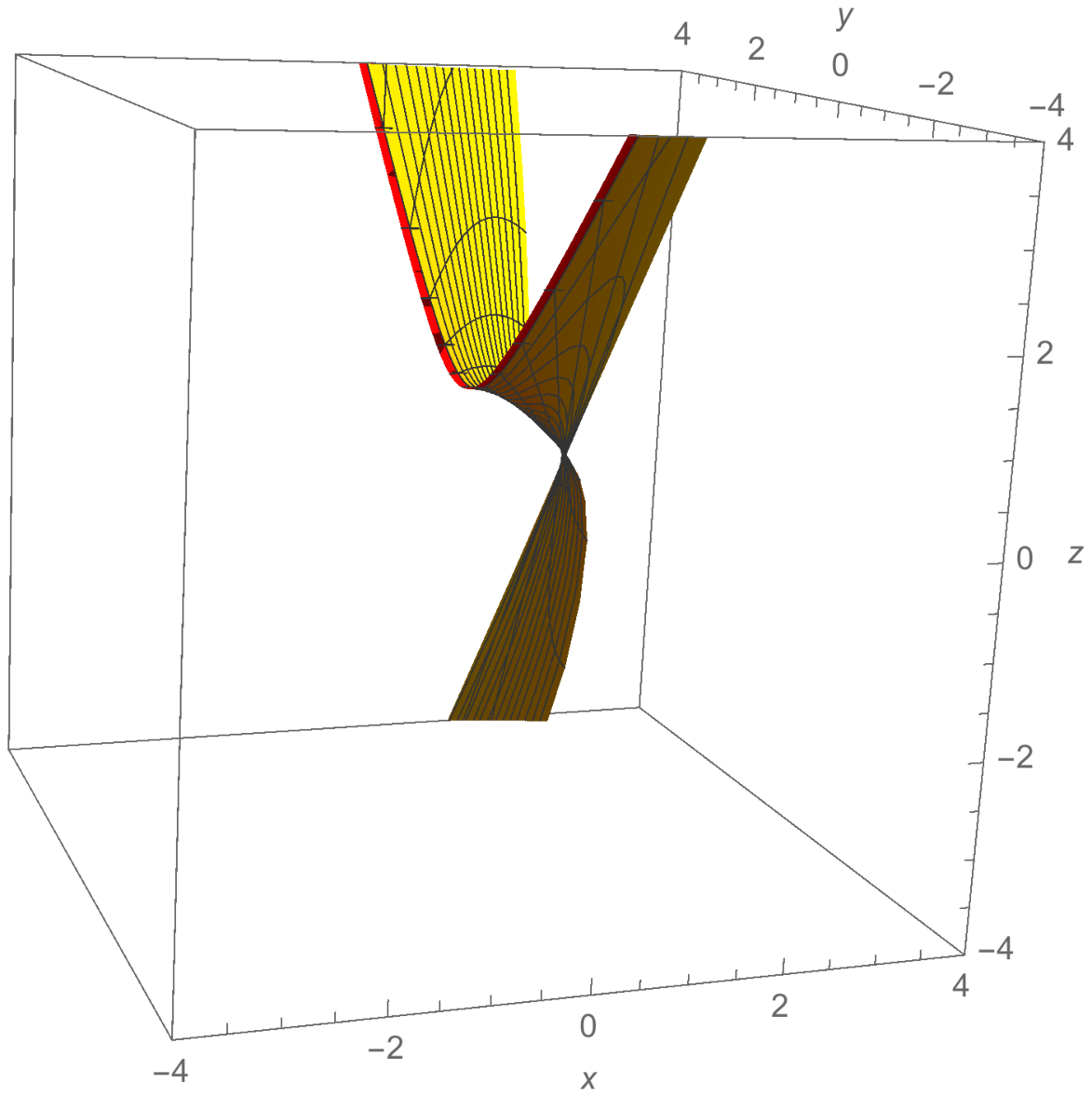}
\includegraphics[width=.4\textwidth]{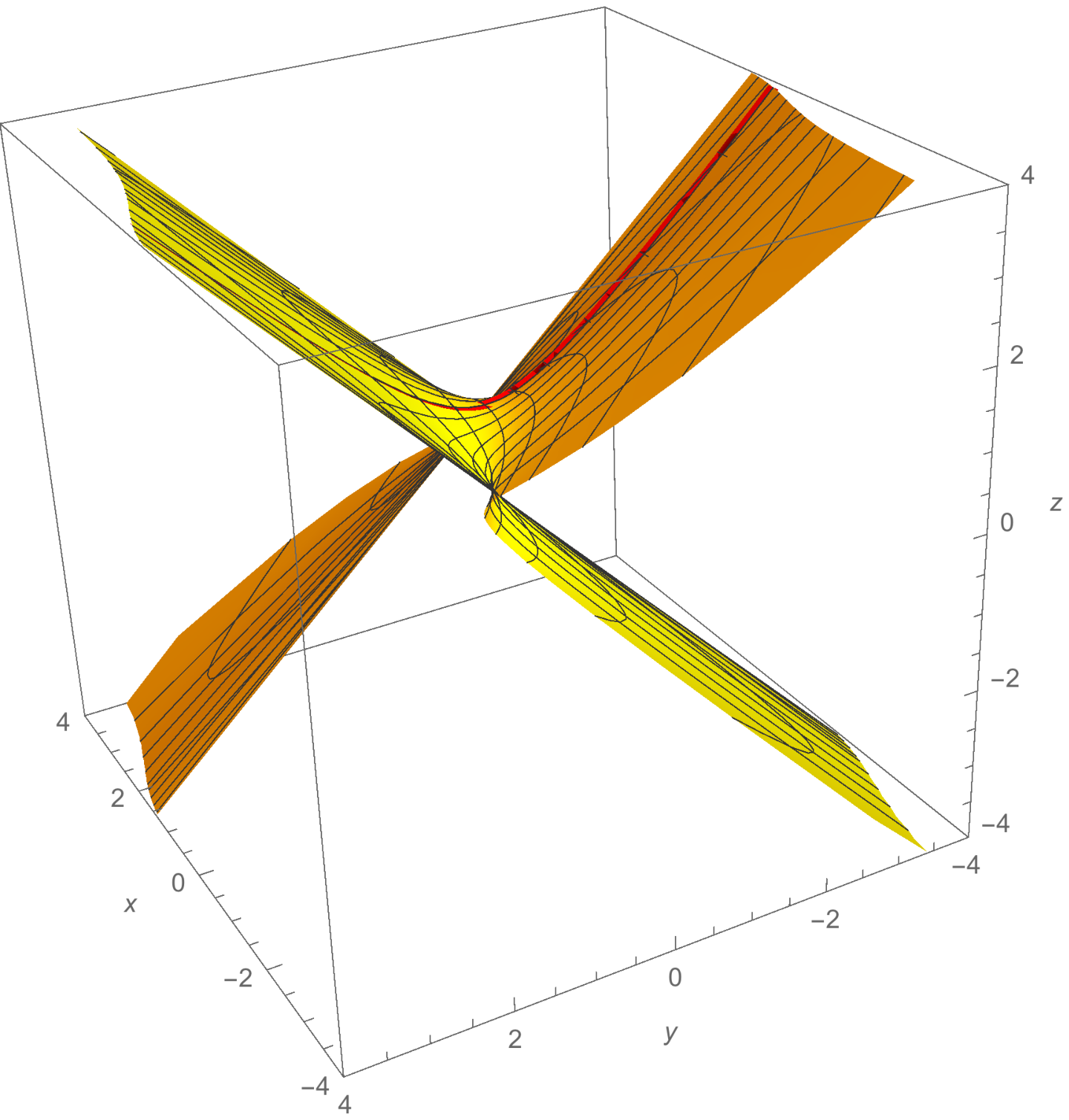}
\caption{A bending helicoid $X(u,v)$ where the core curve is a circle $\alpha$ with spacelike axis. Here $a=1$. Left: a strip of $X(u,v)$ along $\alpha(u)=X(u,0)$. Right: a view of the full surface.}
\end{figure}

\subsection{The axis is  lightlike}\label{sl}
 The parametrization of a   circle with lightlike axis is given in terms of a polynomial function on the parameter $t$ of $\alpha$. Indeed, we can assume that the axis $L$ is spanned by $(1,0,1)$. After an homothety, a circle with axis $L$   parametrizes as $\alpha(t)=(-1+t^2/2,t,t^2/2)$. Then $\alpha'(t)=(t,1,t)$ and $\alpha''(t)=(1,0,1)$. Now $\alpha''(t)$ is a lightlike vector for all $t\in \r$. Under this setting, there is not assigned a Frenet frame of $\alpha$ formed by an orthonormal basis of $\l^3$. Instead, one can define $\n(t)=\alpha''(t)/2$ and $\b(t)$ the unique lightlike vector orthogonal to $\alpha'(t)$ such that $\langle\n(t),\b(t)\rangle=-1/2$ (\cite{lop}). In our case, 
\begin{equation}\label{lig0}
\n(t)=\left(\frac12,0,\frac12\right),\quad\quad\b(t)=\left(\frac{ t^2-1}{2},t,\frac{t^2+1}{2}\right).
\end{equation}
 Then $e_2(t)=\n(t)-\b(t)$ and $e_3(t)=\n(t)+\b(t)$, together $\alpha'(t)$ form an orthonormal basis $\{\alpha'(t),e_2(t),e_3(t)\}$ where $e_3(t)$ is a timelike vector, namely,  
  $$e_2(t)=\left(\frac{2-t^2}{2},-t,-\frac{t^2}{2}\right),\quad\quad e_3(t)=\left(\frac{t^2}{2},t,\frac{t^2+2}{2}\right).$$
The vector field $V(t)$ in the   Bj\"{o}rling  problem writes as 
 $V(t)=\sinh(\varphi(t)) e_2(t)+\cosh(\varphi(t)) e_3(t)$.  There is an important difference with the above cases of timelike and spacelike axis when  we take $\varphi(t)=at$. Now the integrand in \eqref{bj} is formed by trigonometric functions and polynomials on the variable $t$ which it does not allow to solve explicitly the integral of \eqref{bj} by algebraic operations. In this subsection we only consider the case that $\varphi$ is a constant function, namely,   $\varphi(t)=a$, $a\in\r$. Then$$V(z)\times\alpha'(z)=\left(\begin{array}{l}
\frac{1}{2} \left(\cosh (a)(z^2-2)- \sinh (a) z^2\right)\\
 (\cosh (a)-\sinh (a))z\\
\frac{1}{2} \left( \cosh (a)z^2-\sinh (a)(z^2+2) \right)\end{array}\right).$$
The integration of \eqref{bj} gives
$$X(z)=\mbox{Re} \left(\begin{array}{l}
 \frac{1}{6}i z^3 (\cosh(a)-\sinh (a)) -iz \cosh (a) +\frac{z^2}{2}-1\\
z+i \left(\frac{1}{2} z^2 \cosh (a)-\frac{1}{2} z^2 \sinh (a)\right)\\
\frac{z^2}{2}  -\frac{1}{6} iz^3 (\sinh (a)- \cosh (a))-iz \sinh (a)\end{array}\right).$$
Doing $z=u+iv$, we obtain
\begin{equation}\label{lig2}
X(u,v)=\left(
\begin{array}{l} 
(\sinh (a)-  \cosh (a))( \frac{1}{2} u^2 v -\frac{1}{6} v^3 )+v \cosh (a)+\frac{u^2}{2}-\frac{v^2}{2}-1\\
  u+ (\sinh (a)- \cosh (a))uv\\
(\sinh (a)-  \cosh (a))(\frac{1}{2} u^2 v  -\frac{1}{6} v^3 )+v \sinh (a)+\frac{u^2}{2}-\frac{v^2}{2}
\end{array}\right).
\end{equation}

We prove that this surface is a     rotational with respect to the  axis $L$   spanned by $(1,0,1)$. The uniparametric group of rigid motions with  axis $L$  is
$$G_l=\left\{\Psi^l(\theta)=
\left(
\begin{array}{ccc}
 1-\frac{\theta^2}{2} & \theta & \frac{\theta^2}{2} \\
 -\theta & 1 & \theta \\
 -\frac{\theta^2}{2} & \theta & \frac{\theta^2}{2}+1 \\
\end{array}
\right):\theta\in\r\right\}.$$
Each orbit of a point under $\Psi^l(\theta)$   meets the $xz$-plane, so a rotational surface with axis $L$ parametrizes as 
$Z(u,s)=\Psi^l(u)\cdot \beta(s)$, where $\beta$ is curve in the $xz$-plane. Consider $\beta(s)=(s,0,f(s))$. The zero mean curvature condition on $Z(u,s)$ writes as
$(s-f)f''=(f'^2-1)(f'-1)$. An integration of $H=0$ gives  
\begin{equation}\label{lig1}
c(s-f)^3+s-f=2s+b,\quad c>0, b\in\r.
\end{equation}
If we do a change of variables in the generating curve $\beta$ by 
$$\beta(s)=(s,0,f(s))= (h(v)+v,0,h(v)-v),$$ 
then the solution of \eqref{lig1} is $h(v)=\lambda v^3+\mu$, where $\lambda,\mu\in\r$ and $\lambda>0$. This surface is the Enneper's surface of second kind  (see \cite{ko}). Take $\lambda=8\mu+4$ and 
$$\mu=-\frac{1}{3} (\cosh (a)-2 \sinh (a)) (\sinh (a)+\cosh (a)).$$
Then it straightforward to see that the surface $Z(u,v)=\Psi^l(u)\cdot\beta(v)$, where 
$\beta(v)=(v+\lambda v^3+\mu,0,-v+\lambda v^3+\mu)$ is a maximal surface with 
$Z(u,-\frac12)=\alpha(u)$ and the unit normal vector to the surface $Z(u,v)$ along $\alpha$ ($v=-\frac12$) is $V(u)$. By uniqueness of solutions of the Bj\"{o}rling  problem, the surface $Z(u,v)$ must be $X(u,v)$, proving that the surface $X(u,v)$ is a rotational surface.

If we denote $M_{a,l}$ the   surface whose parametrization is \eqref{lig2}, the holomorphic $1$-forms of $M_{a,l}$ are   
$$\phi=\left(\begin{array}{l}
 z+\frac{1}{2} i \left(\left(z^2-2\right) \cosh (a)-z^2 \sinh (a)\right)\\
 1+i z (\cosh (a)-\sinh (a))\\
 z+\frac{1}{2} i \left(z^2 \cosh (a)-\left(z^2+2\right) \sinh (a)\right)
 \end{array} \right)dz.$$

We now consider the minimal surface $M_{a,l}^\sharp\subset\e^3$.  Then  $\psi_k$ have    no real periods and its  Weierstrass representation is
 $$\omega=-\frac{1}{2}  \left( \cosh \left(\frac{a}{2}\right)(z-2i)-\sinh \left(\frac{a}{2}\right)z \right)^2dz,$$
 $$g(z)=\frac{ \sinh \left(\frac{a}{2}\right)(iz-2)- i z \cosh \left(\frac{a}{2}\right)}{  \cosh \left(\frac{a}{2}\right)(z-2i)-\sinh \left(\frac{a}{2}\right)z}.$$
This surface has only one end, namely, $z=\infty$, and the order of the pole at $\infty$ is $4$. It is also clear that $M_{a,l}$   is a complete   with total curvature $-4\pi$ and by the classification of complete minimal surfaces of total curvature $-4\pi$ (\cite{os}), this  surface is the Enneper's surface.

\section{Bj\"{o}rling surfaces based on a helix }
We consider the Bj\"{o}rling problem when the basis curve $\alpha$ is a spacelike helix of $\l^3$. We only consider that the axis of the helix is timelike or spacelike.  By a  \emph{helix}   in $\l^3$ we mean the orbit of a point under a uniparametric group of helicoidal  motions  when this orbit is not a straight line. Again, we will take $V(t)$ as a  linear combination of the normal and binormal vectors formed by trigononometric polynomials. In Euclidean setting, these surfaces are called helicoidal helicoids \cite{we}.

\subsection{Timelike axis}

Up to a homothety, a spacelike helix with   axis $(0,0,1)$ parametrizes as $\alpha(t)=(\cos(t),\sin(t),\lambda t)$, with $0<\lambda<1$. Although $t$ is not the arc-length parameter of   $\alpha$, the normal vector $\n(t)$ and binormal vector $\b(t)$ are
$$\n(t)=(-\cos(t),-\sin(t),0),\quad\quad\b(t)=\frac{1}{\mu}(\lambda\sin(t),-\lambda\cos(t),-1),$$
with $\mu=\sqrt{1-\lambda^2}$. A unit timelike vector field  $V(t)$ expresses as $V(t)=\sinh(\varphi(t)) \n(t)+\cosh(\varphi(t)) \b(t)$. The case that $\varphi$ is a constant function gives a helicoidal surface. Indeed, if $\varphi(t)=a\in\r$, then an integration of \eqref{bj} gives
$$X(u,v)=\left(\begin{array}{l}
-\mu \cosh (a) \cos (u) \sinh (v)+\lambda  \sinh (a) \sin (u) \sinh (v)+\cos (u) \cosh (v)\\
\cosh(v)\sin(u)-\mu \cosh (a) \sin (u) \sinh (v)-\lambda  \sinh (a) \cos (u) \sinh (v)\\
\lambda  u-v \sinh (a)\end{array}
\right).$$
For any $\theta\in\r$, consider the helicoidal motion of axis $(0,0,1)$
$$\Phi^t(\theta):(x,y,z)\longmapsto \left(
\begin{array}{ccc}
\cos(\theta)  & -\sin(\theta)&0\\
\sin(\theta) & \cos(\theta) & 0 \\
0 & 0&1 \\
\end{array}
\right)\left(\begin{array}{l}x\\ y\\  z\end{array}\right)+  \left(\begin{array}{l}0\\ 0\\  \lambda \theta\end{array}\right).$$
Then it is immediate that $\Phi^t(\theta)\cdot X(u,v)=X(u+\theta,v)$ for all $\theta\in\r$, proving that the surface is helicoidal

Suppose now that $\varphi(t)=at$, $a>0$. Then
$$V(t)=\left(\begin{array}{c}\frac{\lambda}{ \mu} \sin (t) \cosh (a t)-\cos (t) \sinh (a t)\\
- \frac{\lambda}{ \mu} \cos (t) \cosh (a t) -\sin (t) \sinh (a t)\\
-\frac{1}{ \mu}\cosh (a t),\end{array}\right),\quad\quad t\in\r.$$

By integrating   \eqref{bj} we obtain the expression of $X(u,v)$, namely, 
\begin{equation}
\begin{aligned}\label{ht}
&\Big( 
 \frac{\cosh (v) \left(\left(a^2+1\right) \cos (u)-\left(\mu -a \lambda \right) \sin (u) \sinh (a u) \sin (a v)-\left(a \mu +\lambda \right) \cos (u) \cosh (a u) \sin (a v)\right)}{a^2+1}\\
 &+\frac{\sinh (v) \cos (a v) \left(\left(a \mu +\lambda \right) \sin (u) \sinh (a u)-\left(\mu -a \lambda \right) \cos (u) \cosh (a u)\right)}{a^2+1},\\
 &\frac{\cosh (v) \left(\left(a^2+1\right) \sin (u)-\left(a \mu +\lambda \right) \sin (u) \cosh (a u) \sin (a v)+\left(\mu -a \lambda \right) \cos (u) \sinh (a u) \sin (a v)\right)}{a^2+1}\\
 &-\frac{\sinh (v) \cos (a v) \left(\left(a \mu +\lambda \right) \cos (u) \sinh (a u)+\left(\mu -a \lambda \right) \sin (u) \cosh (a u)\right)}{a^2+1},\\
 &\lambda  u-\frac{\sinh (a u) \sin (a v)}{a} \Big)
\end{aligned}
\end{equation}
We now give the Weierstrass data of $X(u,v)$. The holomorphic $1$-forms $\phi_k$ are
\begin{eqnarray*}
\phi_1&=&\left( i \mu  \cos (z) \cosh (a z)-\sin (z) (1+i \lambda  \sinh (a z))\right) dz\\
\phi_2&=&\left(\cos(z)(1+i\lambda\sinh(az))+i \mu  \sin (z) \cosh (a z)\right) dz\\
\phi_3&=&\left(\lambda+i\sinh (az)\right)dz,
\end{eqnarray*}
and the Weierstrass data are
$$\omega=\frac{1}{2} e^{-(a+i) z} \left(\left(\lambda +i \mu\right) e^{2 a z}-2 i e^{a z}+i \mu-\lambda \right)dz.$$
$$g(z)=\frac{e^{i z} \left(-2 i \lambda  e^{a z}+e^{2 a z}-1\right)}{\left(\mu-i \lambda \right) e^{2 a z}-2 e^{a z}+\mu+i \lambda }.$$

\begin{figure}[h]
\includegraphics[width=.4\textwidth]{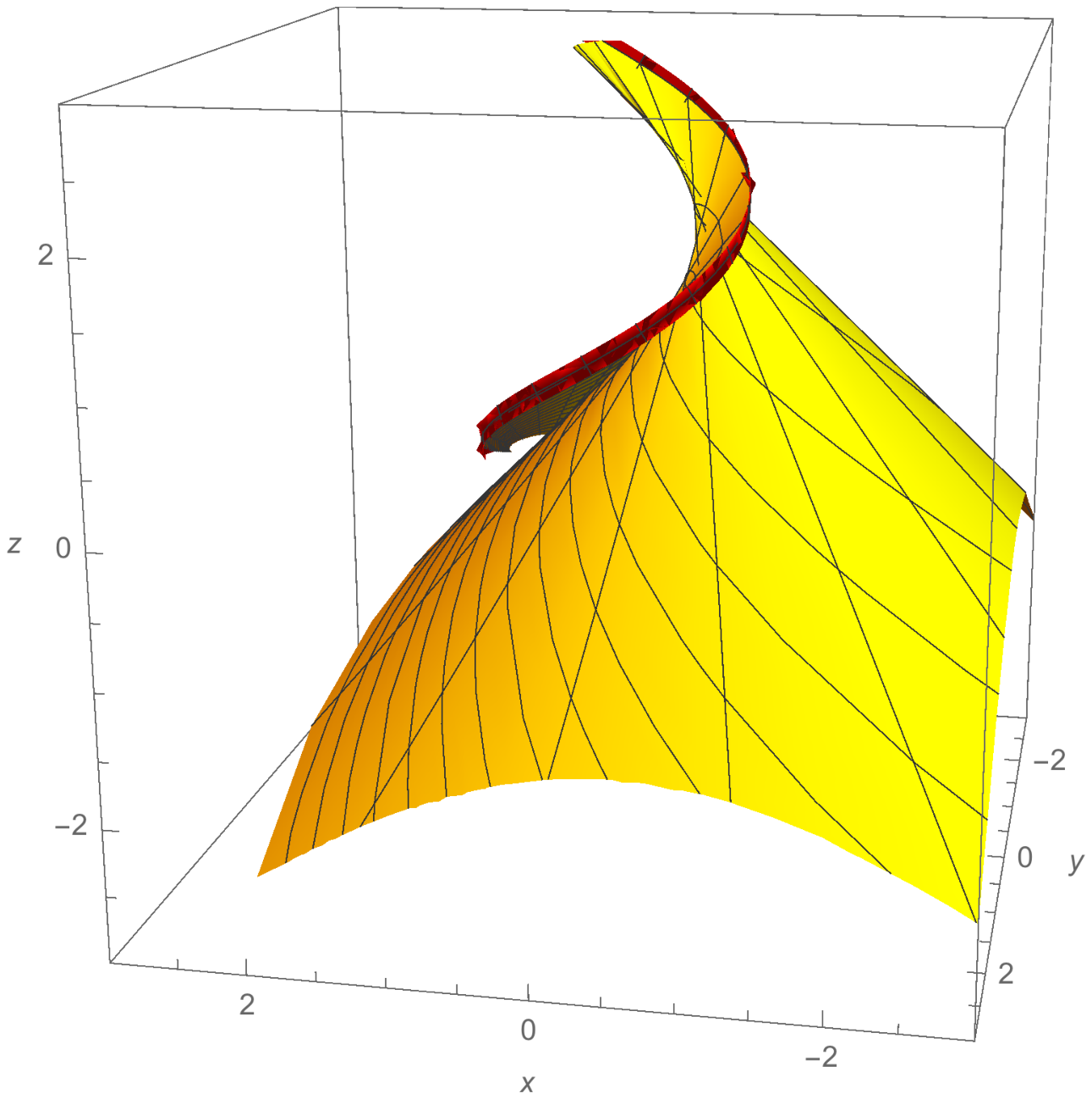}
\includegraphics[width=.4\textwidth]{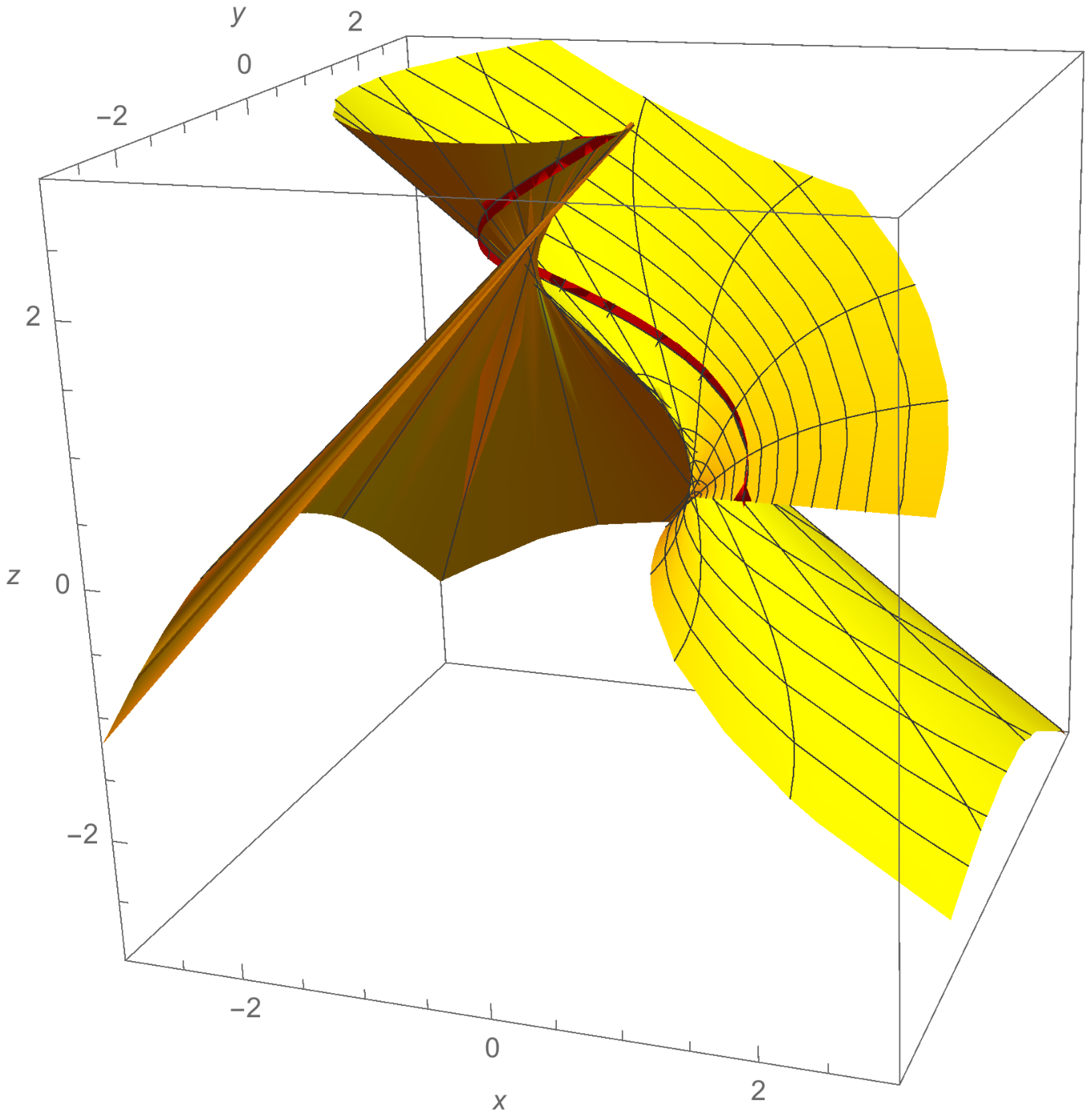}
\caption{A helicoidal helicoid $X(u,v)$ where the core curve is a circle $\alpha$ with timelike axis. Here $a=1$, $\lambda=3/5$ and $\mu=4/5$. Left: a strip of $X(u,v)$ along $\alpha(u)=X(u,0)$. Right: a view of the full surface.}
\end{figure}
\subsection{Spacelike axis}

  Suppose that the axis $L$ is $(1,0,0)$. Then there are two types of spacelike helices:  
\begin{eqnarray*}
\alpha(t)&=&(\lambda t,\cosh(t),\sinh(t)),\quad \lambda>1, \mbox{(type I)}\\
\alpha(t)&=&(\lambda t,\sinh(t),\cosh(t)), \quad \lambda>0, \quad \mbox{(type II)}.
\end{eqnarray*}
Because the computations are similar as in the case of timelike axis, we only give the parametrizations of the surface and its Weierstrass representation. Again we point out that the curve $\alpha$ is not parametrized by the arc-length.  As in the case of timelike axis, if the function $\varphi$ that defines the unit timelike vector field $V(t)$ is a constant function, the Bj\"{o}rling surface   is a helicoidal surface. We omit the details. 
\subsubsection{Helix of type I}

 The normal vector $\n(t)$ and binormal vector $\b(t)$ are
$$\n(t)=(0,\cosh(t),\sinh(t)),\quad\quad\b(t)=-\frac{1}{\mu}(1,\lambda\sinh(t),\lambda\cosh(t)),$$
where $\mu=\sqrt{\lambda^2-1}$. Because $\n(t)$ is spacelike and $\b(t)$ is timelike, the unit timelike vector field  $V$   writes as $V(t)=\sinh(\varphi(t))\n(t)+\cosh(\varphi(t))\b(t)$. Consider the case that $\varphi(t)=at$. When $a\not=1$, the Bj\"{o}rling surface $X(u,v)$ is
\begin{equation}
\begin{aligned}\label{hs11}&\Big( 
 \lambda  u-\frac{\sinh (a u) \sin (a v)}{a},\\
& \frac{\cosh (u) \left(\left(a^2-1\right) \cos (v)+\cosh (a u) \left(\left(a \mu +\lambda \right) \cos (v) \sin (a v)-\left(a \lambda +\mu \right) \sin (v) \cos (a v)\right)\right)}{a^2-1}\\
&+\frac{\sinh (u) \sinh (a u) \left(\left(a \mu +\lambda \right) \sin (v) \cos (a v)-\left(a \lambda +\mu \right) \cos (v) \sin (a v)\right)}{a^2-1},\\
 &\frac{\sinh (u) \left(\left(a^2-1\right) \cos (v)+\cosh (a u) \left(\left(a \mu +\lambda \right) \cos (v) \sin (a v)-\left(a \lambda +\mu \right) \sin (v) \cos (a v)\right)\right)}{a^2-1}\\
 &+\frac{\cosh (u) \sinh (a u) \left(\left(a \mu +\lambda \right) \sin (v) \cos (a v)-\left(a \lambda +\mu \right) \cos (v) \sin (a v)\right)}{a^2-1} 
 \Big)
\end{aligned}
\end{equation}
and if $a=1$,   the Bj\"{o}rling surface  $X(u,v)$ is
\begin{equation}
\begin{aligned}\label{hs12}&\Big( 
 \lambda  u-\sinh (u) \sin (v),\\
& \frac{1}{2} \left(\left(\mu -\lambda \right) (\sinh ^2(u)+\cosh^2(u)) \sin (v) \cos (v)+2 \cosh (u) \cos (v)+\left(\mu +\lambda \right) v\right),\\
&\sinh (u) \cos (v) \left(\left(\mu -\lambda \right) \cosh (u) \sin (v)+1\right) \Big)
\end{aligned}
\end{equation}

The holomorphic functions $\phi_k$ are
\begin{eqnarray*}
\phi_1&=&\left(\lambda +i \sinh (a z)\right) dz\\
\phi_2&=&\left(\sinh (z)+i(\lambda  \sinh (z) \sinh (a z) - \mu  \cosh (z) \cosh (a z))\right)dz\\
\phi_3&=&\left(\cosh (z)+i(\lambda  \cosh (z) \sinh (a z)- \mu  \sinh (z) \cosh (a z))\right) dz.
\end{eqnarray*}
We write the above trigonometric functions in terms of the exponential function  and by the change of variable $e^z\rightarrow z$ we have
\begin{eqnarray*}
\phi_1&=&\frac{i z^{2 a}+2 \lambda  z^a-i}{2 z^{a+1}}dz\\
\phi_2&=&\frac{i (\lambda-\mu)  (z^{2 a+2}+1)-i (\lambda+\mu)  (z^{2 a}+z^2)+2 z^{a+2}-2 z^a  }{4 z^{a+2}  }dz\\
\phi_3&=& \frac{i (\lambda-\mu)  (z^{2 a+2}-1)+i (\lambda+\mu)  (z^{2 a}-z^2)+2 z^{a+2}+2 z^a  }{4 z^{a+2}  }dz
\end{eqnarray*}
Suppose $a=n\in\mathbb{N}$. Then it is clear that if $n\not=1$, the $1$-forms $\phi_k$ have no periods and if $n=1$, then only $\phi_2$ has real periods. The Weierstrass data are
$$\omega=\frac{ (\lambda-\mu)(z^{2n+2}+1)-\left(\mu+\lambda \right)( z^{2 n}+z^2)+2i (z^{n+1}-z^{n+2}+z^n-z) +4 \lambda  z^{n+1}  }{4 z^{n+2}}dz.$$
$$g(z)=\frac{-\left(\mu+\lambda \right) (z^{2 n}-z^2)+\left(\mu-\lambda \right) (z^{2 n+2}-1)+2 i (z^{n+2}+  z^n ) }{-i \left(\mu+\lambda \right) (z^{2 n}+z^2)-i \left(\mu-\lambda \right) (z^{2 n+2}+1)+4 i \lambda  z^{n+1} -2 (z^{2 n+1}-z^{n+2}+z^n-z)  }.$$

\subsubsection{Helix of type II} 
Now   the normal vector $\n(t)$ and binormal vector $\b(t)$ are given by
$$\n(t)=(0,\sinh(t),\cosh(t)),\quad\quad\b(t)=\frac{1}{\mu}(1,-\lambda\cosh(t),-\lambda\sinh(t)),$$
where $\mu=\sqrt{\lambda^2+1}$. A unit timelike vector field  $V(t)$ orthogonal to $\alpha$ writes as $V(t)=\cosh(\varphi(t))\n(t)+\sinh(\varphi(t))\b(t)$. Consider the case that $\varphi(t)=at$.  When $a\not=1$,  the Bj\"{o}rling surface $X(u,v)$ is 
\begin{equation}
\begin{aligned}\label{hs21}&\Big( 
\frac{\cosh (a u) \sin (a v)}{a}+\lambda  u,\\
& \frac{\sinh (u) \left(\left(a^2-1\right) \cos (v)+\sinh (a u) \left(\left(a \mu +\lambda \right) \cos (v) \sin (a v)-\left(a \lambda +\mu \right) \sin (v) \cos (a v)\right)\right)}{a^2-1}\\
&+\frac{\cosh (u) \cosh (a u) \left(\left(a \mu +\lambda \right) \sin (v) \cos (a v)-\left(a \lambda +\mu \right) \cos (v) \sin (a v)\right)}{a^2-1},\\
&\frac{\cosh (u) \left(\left(a^2-1\right) \cos (v)+\sinh (a u) \left(\left(a \mu +\lambda \right) \cos (v) \sin (a v)-\left(a \lambda +\mu \right) \sin (v) \cos (a v)\right)\right)}{a^2-1}\\
&+\frac{\sinh (u) \cosh (a u) \left(\left(a \mu +\lambda \right) \sin (v) \cos (a v)-\left(a \lambda +\mu \right) \cos (v) \sin (a v)\right)}{a^2-1} 
 \Big)\end{aligned}
 \end{equation}
 and if $a=1$, the surface is 
\begin{equation}
\begin{aligned}\label{hs22}&\Big( 
\lambda  u+\cosh (u) \sin (v),\\
& \frac{1}{2} \left(\left(\mu -\lambda \right) (\sinh ^2(u)+\cosh^2(u)) \sin (v) \cos (v) +2 \sinh (u) \cos (v)-\left(\mu +\lambda \right) v\right),\\
&\cosh (u) \cos (v) \left(\left(\mu -\lambda \right) \sinh (u) \sin (v)+1\right).
 \Big)\end{aligned}
 \end{equation}
 
The holomorphic functions $\phi_k$ are
\begin{eqnarray*}
\phi_1&=&\left(\lambda -i \cosh (a z)\right)dz\\
\phi_2&=&\left(\cosh (z)+i( \lambda  \cosh (z) \cosh (a z)- \mu \sinh (z) \sinh (a z))\right) dz\\
\phi_3&=&\left(\sinh (z)+i( \lambda  \sinh (z) \cosh (a z)- \mu \cosh (z) \sinh (a z))\right)dz.
\end{eqnarray*}
Again we do the change of variables $e^{iz}\rightarrow z$, obtaining
\begin{eqnarray*}
\phi_1&=&\frac{-i z^{2 a}+2 \lambda  z^a-i}{2z^{a+1}}dz\\
\phi_2&=&\frac{i (\lambda-\mu)  (z^{2 a+2}+1)+i (\lambda+\mu) (z^{2 a}+z^2)+2 z^{a+2}+2 z^a  }{4z^{a+2}}dz \\
\phi_3&=&\frac{i (\lambda-\mu) (z^{2 a+2}-1)-i (\lambda+\mu)  (z^{2 a}-z^2)+2 z^{a+2}-2 z^a }{4z^{a+2}}dz.\end{eqnarray*}
 Consider $a=n\in\mathbb{N}$. Again,   the holomorphic $1$-forms $\phi_k$  have no real periods if $n\not=1$. When $n=1$, then only $\phi_2$ has real periods.  The Weierstrass data are
 $$\omega=\frac{\left(\lambda -\mu\right) (z^{2 n+2}+1)+  \left(\mu+\lambda \right) (z^{2 n}+z^2)+4 \lambda  z^{n+1}-2 i (z^{n+2} +  z^{2 n+1} +  z^n  + z)}{4 z^{n+2}}dz,$$
 $$g(z)=\frac{i \left(\mu+\lambda \right) (z^{2 n}-z^2)+i \left(\mu-\lambda \right) (z^{2 n+2}-1)-2 (z^{n+2}- z^n) }{\left(\mu-\lambda \right)( z^{2 n+2}+1)-\left(\mu+\lambda \right) (z^{2 n}+z^2)-4 \lambda  z^{n+1}+2 i (z^{n+2}+  z^{2 n+1}+ z^n +  z)}.$$ 
 
 \begin{figure}[h]
\includegraphics[width=.4\textwidth]{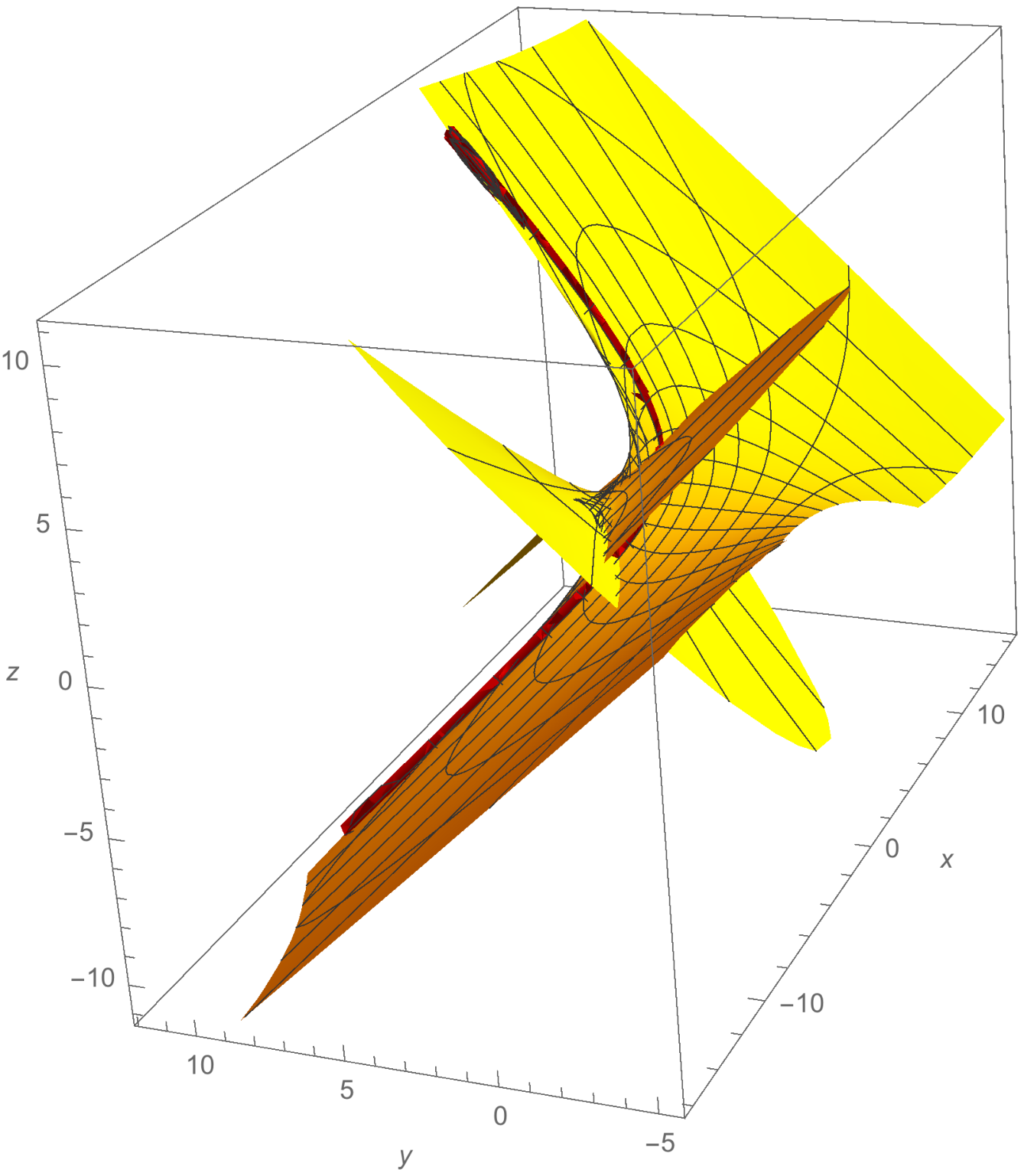}
\includegraphics[width=.4\textwidth]{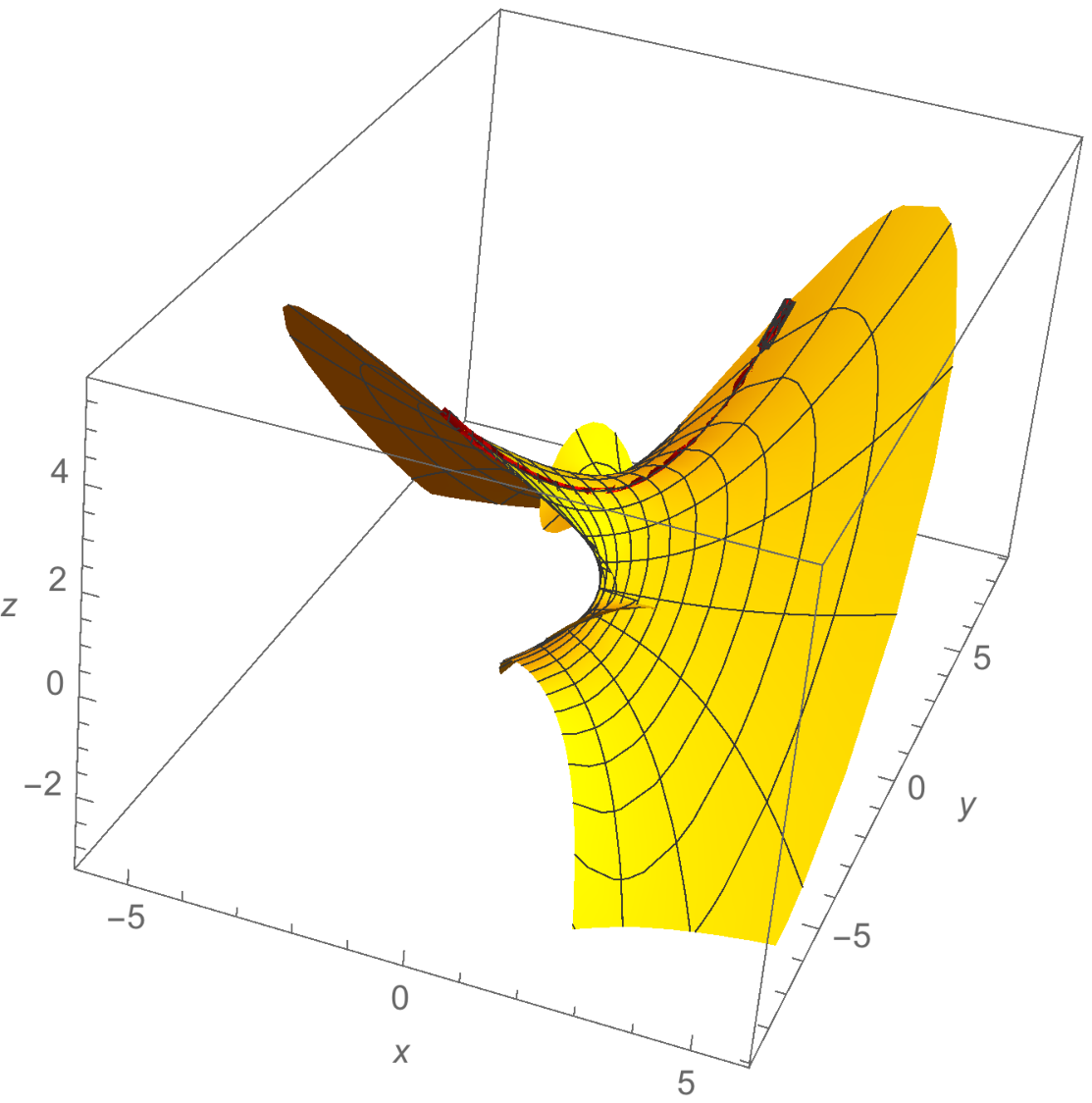}
\caption{A helicoidal helicoid $X(u,v)$ where the core curve is a helix $\alpha$ with spacelike axis.  Left: the helix $\alpha$ is of type I, with $a=1$, $\lambda=2$ and $\mu=\sqrt{3}$.  Right: the helix $\alpha$ is of type II, with $a=1$, $\lambda=1$ and $\mu=\sqrt{2}$.}
\end{figure}

\end{document}